\documentclass[reviewcopy]{elsart}

\usepackage{amssymb}
 \usepackage{epsfig}
\theoremstyle{plain}
\newtheorem{theorem}{Theorem}[section]
\newtheorem{corollary}[theorem]{Corollary}

\newtheorem{proposition}[theorem]{Proposition}
\theoremstyle{definition}
\newtheorem{definition}[theorem]{Definition}
\newtheorem{example}[theorem]{Example}

\newcounter{abc}

\newcounter{ABC}

\newcommand{\N}{{\mathbb N}}
\newcommand{\Z}{{\mathbb Z}}

\begin{document}

\begin{frontmatter}


\title {The reverse engineering problem with probabilities and sequential behavior: Probabilistic
Sequential Networks}


\author{Maria A. Avi\~no-Diaz}
\address{Department of Mathematic-Physics\\ University of Puerto Rico\\
Cayey, Puerto Rico 00736 }
\begin{abstract} The reverse engineering problem with probabilities and sequential behavior
is introducing here, using the expression of an algorithm. The solution is partially founded,
because we solve the problem only if we have a Probabilistic Sequential Network.
 Therefore the  probabilistic structure on sequential dynamical systems is introduced here, the new model
will be called Probabilistic Sequential Network, PSN. The morphisms
of Probabilistic Sequential Networks are defined using two algebraic
conditions, whose imply that the distribution of probabilities in
the systems are close. It is proved here that two homomorphic
Probabilistic Sequential Networks have the same equilibrium or
steady state probabilities. Additionally, the proof of the set of
PSN with its morphisms form the category \textbf{ PSN}, having the
category of sequential dynamical systems \textbf{ SDS}, as a full
subcategory is given. Several examples of morphisms, subsystems and
simulations are given.
\end{abstract}

\end{frontmatter}

\section{Introduction}\label{intro}

  Probabilistic Boolean Networks  was introduced by I. Schmulevich,
 E. Dougherty, and W. Zhang in 2000, for studying the dynamic of a network using
  time discrete Markov chains, see \cite{S1, S2,SDZ, S3}. This model had  several applications
  in the study of cancer, see \cite{S4}. It is important for development  an
  algebraic mathematical theory of the model Probabilistic Boolean Network
  PBN, to describe special maps between two PBN,
  called homomorphism and projection, the first papers in this direction are, \cite{DS,ID},
  We will use the acronym PBN, PSN, or SDS for plural as well as
  singular instances.
   Instead of this model is being used in applications,
  the connection  of these two digraphs of the model: the graph of genes and
  the State Space is an interesting problem to study. The introduction of probabilities
  in the definition of  Sequential Dynamical System has this
  objective. This paper is  the first  part of this theory.

  The theory of sequential dynamical systems (SDS) was born studying
  networks where the entities involved in the problem do not necessarily arrive
  at a place  at the same time, and it is part of the theory of computer
  simulation, \cite{B1, B2}. The mathematical background for SDS was recently development  by
Laubenbacher and Pareigis, and  it solves  aspects of the theory and
applications, see \cite{LP1,LP2, LP3}.

The introduction of a probabilistic structure on Sequential
Dynamical Systems is an  interesting problem that it is introduced
in this paper. A SDS induces a finite dynamical system $(k^n,f)$,
for the classifications of Linear Dynamical Systems see \cite{H},
but the mean difference between a SDS and FDS
  is that there exits another graph with new information giving by the
  local functions, and an order in the sequential behavior of these
  local functions. It is known, that a finite dynamical
  systems can be studied as a SDS, because we can construct a bigger
  system that in this case is sequential.
  Making together the sequential order and the probabilistic structure in the dynamic of the system,
   the possibility to work in applications to genetics increase, because genes act in a
  sequential manner. In particular the notion of morphism in the category of SDS
  establishes connection between the digraph of genes and the State Space, that is the dynamic of the function.
   Working in the applications,  Professor Dougherty's group  wanted to consider two things in the definition of PBN:
  a sequential behavior on genes, and the exact definition of projective maps
  between two PBN that inherits the properties of the first digraph of genes.
   For this reason, a new model that considers both questions and  tries to construct
  projections that  work well is described here. I introduce in this paper
  the sequential behavior and the probability together in PSN and
  my final objective is to construct  projective maps that let us
  reduce the number of functions in the finite dynamical systems
  inside the PBN. One of the mean problem in modeling dynamical systems
  is the computational aspect of the number of functions and the
  computation  of steady states in the State Space. In particular, the reduction of number of functions is one of
the most important problems, because by solving that we can
determine which part of the network \emph{State Space} may be
simplified.  The concept of morphism, simulation, epimorphism, and
equivalent Probabilistic Sequential Networks are developed  in this
paper, with this particular objective.

This paper is organized  as follows. In section \ref{SDS}, a
notation slightly different to the one used in \cite{LP2}  is
introduced for homomorphisms of SDS. This notation is helpful for
giving the concept of morphism of PSN. In section \ref{PSN1}, the
probabilistic structure on SDS is introduced using for each vertex
of the support graph, a set of local functions,  more than one
schedule, and finally having several update functions with
probabilities assigned to them. So, it is obtained a new concept:
probabilistic sequential network (PSN). The concept of morphism of
PSN is defined with two conditions, one of the most interesting
results in this paper is that these algebraic conditions implies  a
probabilistic condition about the distribution of probabilities in
both PNS, it is proved in Theorem \ref{C3}. In Theorem \ref{MT} is
proved that  two homomorphic PSN have the same \emph{equilibrium or
steady state probabilities}. These strong results justify the
introduction of the dynamical model PSN as an application to the
study of sequential systems.  In section \ref{CAT}, we prove that
the PSN with its morphisms form the category \textbf{ PSN}, having
the category \textbf{ SDS} as a full subcategory. Several examples
of morphisms, subsystems and simulations are given in Section
\ref{EHOM}.
\section{Preliminaries}\label{SDS}
In this introductory section we give the  definitions and results of
Sequential Dynamical System  introduced  by Laubenbacher and
Pareigis in \cite{LP2}. Let $\Gamma$ be a graph, and let
$V_\Gamma=\{1,\ldots, n\}$ be the set of vertices of $\Gamma$. Let
$(k_i|i \in V_\Gamma)$ be a family of finite sets. The set $k_a$ are
called the set of local states at $a$, for all $a\in V_\Gamma$.
Define $ k^{n} := k_1 \times \cdots \times k_n$ with
  $|k_i|<\infty $, the set of (global) states of $\Gamma$.\\
\subsection{Sequential Dynamical System}\label{SDS1}
 A Sequential Dynamical System (SDS) ${\mathcal F}=(\Gamma,(k_i)_{i=1}^{n},(f_i)_{i=1}^n, \alpha)$ consists of
\begin{enumerate}
\item [1.]  A finite graph
$\Gamma=(V_\Gamma,E_\Gamma)$ with the set of vertices
$V_\Gamma=\{1,\ldots ,n\}$ , and the set of edges $E_\Gamma
\subseteq V_\Gamma \times V_\Gamma$.
 \item[2.] A family of finite  sets $(k_i|i\in V_\Gamma)$.
 \item[3.]  A family of  local functions $f_i:
k^n \rightarrow k^n $,  that is  \[f_i(x_1, \ldots , x_n)= (x_1,
\ldots ,x_{i-1}, \overline{f}, x_{i+1}, \ldots , x_n)\] where
$\overline{f}(x_1,\ldots , x_n)$ depends only of those variables
which are connected to $i$ in $\Gamma$.
\item [4.] A
permutation $\alpha =( \begin{array}{lll} \alpha_1& \ldots &
\alpha_n
\end{array}) $ in the set of vertices $V_\Gamma$, called an update
schedule ( i.e. a bijective map $\alpha :V_\Gamma\rightarrow
V_\Gamma)$.
\end{enumerate}
\par The global update function of the SDS is $f=f_{\alpha_1}\circ\ldots\circ{f_{\alpha_n}}$.
 The function $f$ defines the dynamical behavior of the SDS and determines a finite directed graph
  with vertex set $k^n$ and directed edges $(x,f(x))$, for all $x\in {k}^n$, called the State
  Space of $\mathcal{F}$, and denoted by $\mathcal S_f $.\\
The definition of homomorphism between two SDS uses the fact that
the vertices $V_\Gamma=\{1,\ldots , n\}$ of a SDS and the states
$k^n$ together with their evaluation map $k^{n}\times V_{\Gamma} \ni
(x,a) \mapsto <x,a>:=x_a\in k_i$, form a contravariant setup, so
that morphisms between such structures should be defined
contravariantly, i.e. by a pair of certain maps $\phi : \Gamma
\rightarrow \Delta ,$  and the induced function $h_\phi :
k^{m}\rightarrow k^{n}$  with the graph $\Delta$ having $m$
vertices. Here we use  a notation slightly different that the one
using in  \cite{LP2}.

Let $F=(\Gamma,(f_i:k^n\rightarrow {k}^n),\alpha)$ and
$G=(\Delta,(g_i:{k}^m\rightarrow k^m),\beta)$ be two SDS. Let
$\phi:\Delta \rightarrow \Gamma$ be a digraph morphism. Let
$(\widehat{\phi}_b:\overline{}k_{\phi(b)}\rightarrow k_b, \forall
b\in \Delta),$
 be a family of maps in the
category of \textbf{Set}. The map $h_\phi$ is an adjoint map,
because is defined as follows: consider the pairing $k^n\times
V_{\Gamma} \ni (x,a) \mapsto <x,a>:=x_a\in k_a;$ and  similarly $
k^m\times V_\Delta \ni (y,b) \mapsto <y,b>:=y_b\in k_b.$ The induced
adjoint map holds
$<h_\phi(x),b>:=\hat{\phi}_b(<x,\phi(b)>)=\hat{\phi}_b(x_{\phi(b)})$.
Then $\phi,$ and $(\widehat{\phi}_b)$ induce the adjoint map
$h_\phi:k^{n}\rightarrow k^{m}$ defined as follows:
\begin{equation}\label{defh}
           h_\phi(x_1,\ldots,x_n)=(\widehat{\phi}_{1}(x_{\phi(1)}),\ldots,\widehat{\phi}_{m}(x_{\phi(m)})).
\end{equation}

Then $h:F\rightarrow  G$ is a homomorphism of SDS if for all sets of
orders $\tau _{\beta}$ associated to $\beta$ in the connected
components of $\Delta $, the map $h_\phi$ holds the following
conditions:
\begin{equation}\label{eq1}
\left( g_{\beta_l}\circ g_{\beta_{l+1}}\circ \cdots \circ
g_{\beta_s} \right)\circ h_\phi=h_\phi\circ f_{\alpha_i},\
\begin{array}{lcl}
\vspace{-.3in}
 k^n &
\overrightarrow {\hspace{.5in}_{f_{\alpha_i}}\hspace{.5in}  }
& k^n \\
\vspace{-.3in}
\mid && \mid \\
\downarrow  h_\phi &  &  \downarrow  h_\phi \\
\vspace{-.3in}
 k^m &\overrightarrow{\hspace{.2in}_{g_{\beta_l}\circ
g_{\beta_{l+1}}\circ \cdots \circ g_{\beta_s}}\hspace{.3in}}& k^m
\\
\end{array}
\end{equation}
 where $\{\beta_l,\beta_{l+1},\ldots,\beta_s\}=\phi
^{-1}(\alpha_i)$. If $\phi^{-1}(\alpha_i)=\emptyset$, then $Id_{k^m}
\circ h_\phi=h_\phi\circ f_{\alpha_i}$, and the commutative diagram
is now the following: \begin{equation}\label{eq2}
\hspace{1.7in}\begin{array}{lcl}
\vspace{-.2in} k^n & \overrightarrow
{\hspace{.3in}_{f_{\alpha_i}}\hspace{.3in}  }
& k^n \\
\vspace{-.3in}
\mid h_\phi&& \mid h_\phi\\
\vspace{-.3in}
\downarrow  &  & \downarrow  \\
\vspace{-.2in}
 k^m &\overrightarrow {\hspace{.3in}_{Id_{k^m}\hspace{.3in}
 }}& k^m \cr
\end{array}
\end{equation}
 For examples and properties see\cite{LP2}. It that paper,
 the authors proved that the above diagrams implies the following one
\begin{equation}\label{eq3}
\hspace{1.7in}\begin{array}{lcl}
\vspace{-.2in}
k^n & \overrightarrow
{\hspace{.3in}_{f=f_{\alpha_1}\circ \cdots \circ
f_{\alpha_n}}\hspace{.3in}}& k^n \\
\vspace{-.3in}
\mid h_\phi&& \mid h_\phi\\
\vspace{-.3in}
\downarrow  &  & \downarrow  \\
\vspace{-.2in}
k^m &\overrightarrow{\hspace{.2in}_{g=g_{\beta_1}\circ \cdots \circ
g_{\beta_m}}\hspace{.2in}}  & k^m \cr
\end{array}\end{equation}
\subsection{Probabilistic  Boolean Networks} \cite{S1,S2,SDZ,S4}
 The model Probabilistic Boolean  Network
$\mathcal{A}=\mathcal{A}(\Gamma,F,C)$ is defined by the following:
\begin{itemize}
\item [(1)] a finite digraph $\Gamma =(V_\Gamma,
E_\Gamma)$ with  $n$ vertices.
\item [(2)] a family $F=\{F_1, F_2,\ldots,F_n\}$ of  ordered sets
$F_i=\{f_{i1},f_{i2},\ldots,f_{il(i)}\}$ of functions
$f_{ij}:\{0,1\}^n\rightarrow\{0,1\}$, for $i=1,\cdots,n$, and
$j=1,\ldots,l(i)$ called predictors,
\item [(3)]and a family $C=\{c_{ij}\}_{i,j}$, of selection probabilities. The selection probability that the
function $f_{ij}$ is used for the vertex $i$ is $c_{ij}$.
\end{itemize}
\par The dynamic of the model Probabilistic Boolean Network is given
 by the vector  functions $\mathbf
{f}_k=(f_{1k_1}, f_{2k_2},\ldots,f_{nk_n}):\{0,1\}^n\rightarrow
\{0,1\}^n$ for $1\le{k_i}\le{l(i)}$, and $f_{ik_i}\in{F_i}$,  acting
as a transition function.
 Each variable  $x_i\in \{0,1\}$ represents the state of the vertex  $i$.
 All functions are updated synchronously. At every time step, one of the  functions
is selected randomly from the set $F_i$ according to a predefined
probability distribution.  The selection probability that the
predictor $f_{ij}$ is
  used to predict gene $i$ is equal to
\[c_{ij}=P\{f_{ik_i}=f_{ij}\}=\sum_{k_i=j}{p\{{\mathbf
{f}=f}_k\}}.\] There are two digraph structures associated with a
Probabilistic Boolean Network: the low-level graph $\Gamma$, and the
high-level graph which consists of the states of the system and the
transitions between states. The state space $S$ of the network
together with the set of network functions, in conjunction with
transitions between the states and network functions, determine a
Markov chain. The random perturbation makes the Markov chain
ergodic, meaning that it has the possibility of reaching any state
from another state and that it possesses a long-run (steady-state)
distribution. As a Genetic Regulatory Network (GRN), evolves in
time, it will eventually enter a fixed state, or a set of states,
through which it will continue to cycle. In the first case the state
is called a singleton or fixed point attractor, whereas, in the
second case it is called a cyclic attractor. The attractors that the
network may enter depend on the initial state. All initial states
that eventually produce a given attractor constitute the basin of
that attractor. The attractors represent the fixed points of the
dynamical system that capture its long-term behavior. The number of
transitions needed to return to a given state in an attractor is
called the cycle length. Attractors may be used to characterize a
cell's phenotype (Kauffman, 1993) \cite{K}. The attractors of a
Probabilistic Genetic Regulatory Network  (PGRN) are the attractors
of its constituent GRN. However, because a PGRN constitutes an
ergodic Markov chain, its steady-state distribution plays a key
role. Depending on the structure of a PGRN, its attractors may
contain most of the steady-state probability mass \cite{A,MD,Z}.
\section{ The reverse engineering problem with probabilities, and sequential behavior}
\label{Algorithm}
Here, we give a method that permit us to build  sequential  systems
with probabilities assigned to its update functions. This  algorithm made possible to understand the concept of Probability Sequential Network
and it is dedicate to Prof. Rene Hernandez-Toledo.
\subsection{Algorithm: The reverse engineering problem with probabilities, and sequential behavior}\label{algorithm1}
Input:

1.  $n=$ number of entities in the network under studying,
for example 100 genes, and the set of values for each entity, that we denote by $k_a$. \\
2. A set of  relations $\{m_{a,b}\}$ taking $1$ if the entity $a$ is related to the entity  $b$, and $0$  otherwise.\\
3. A set of finite families of states in the network which gives the time series data for one, two or more update functions,
 $A_1=\{(a^1_{i,1},\ldots,a^2_{i,n-1}, a^1_{i,n})|1\leq i\leq m_x \},$
  $\ldots$, and  $A_s=\{(a^s_{i,1},\ldots,a^s_{i,n-1}, a^s_{i,n})|1\leq i\leq
  m_s \}$.\\
4. A set of values $C=\{c_1,\ldots, c_s\}$ with $s$ probabilities obtained in some way
by the experiment or by the time series data. That is $c_1+\cdots+c_s=1,$ and  $c_i\in[0,1]$.
\begin{itemize}
\item[(Alm1)] Creation the low level graph $\Gamma$:
\begin{itemize}
\item[1.]  $V_\Gamma=\{1,\ldots,n\}$ is the set of vertices,
 $k_a$  gives a set of values to each vertex $a$,
 \item[2.] $E_\Gamma=\{(a,b)|if \ m_{a,b}=1  \}$
 We  obtain this graph, using the experiment giving by the
specialists for example, see \cite{BIO}
\end{itemize}
\item [(Alm2)] Denoting $k^n=k_1\times \cdots \times k_n$, we construct
 the local functions $f_{ai}: k^n\rightarrow k^n$  using the data
giving by the experiment, associated to that we have the statistics
of the entities and we give the probability to each function using
the activity of the vertex. Finally we have a set of families of
functions that we denote by $F_a=\{f_{ai}:k^n\rightarrow k^n|1\leq
i \leq \ell(i)\},$ where  $\ell (i)$ is the number of local
functions associated to the vertex $a$.
\item[(Alm3)] The possible types of sequential behavior giving by
the experiment, or we suppose two or three orders  of action with
the possibility to determine which one is better for the network, so
we define $\alpha  = ( \alpha_{1}\hspace{.1in} \ldots \hspace{.1in}
\alpha_{n} )$ in the set of vertices $V_\Gamma $.

\item [(ALM4)] We assigns probabilities to each update function.We
have  The number of update functions is the number of different set $A_j$
 with the time series data for the functions.
\item[(Alm5)] We select a subset of functions such that the behavior
of the network  are closer the experiments, using the probabilities
giving in the set $C=\{c_1, \ldots ,c_s\}$, selected by the
experiments.
\item[(Alm6)] We construct the high level digraph with the selected functions by the set $C$
 in (Alm5).
\end{itemize}
Output: ${\mathcal D}=(\Gamma, \{F_a\}_{a=1}^{|\Gamma|=n}, (k_a)_{a=1}^{n}, (\alpha_j)_{j=1}^m, C=\{c_1,\ldots, c_s\})$.
\begin{example}
Input:\\
1.  $n=3$;  $k_1=k_2=k_3=\Z_2=\{0,1\}$.\\
2.  $\{m_{1,2}=0,\ m_{2,3}=1,\ m_{1,3}=0\}$.\\
3. The data for two update functions: $A_1=\{(0,1,0); (1,1,1); (1,1,0);(1,1,1)\},$\\
  $A_2=\{(0,1,1);(0,1,0); (0,1,1);(0,1,0)\}.$ \\
4. We assign the following probabilities to the function
$C=\{2/3,1/3\} $, using statistic see for example \cite{COE},
\cite{S1}.

Running the algorithm, we obtain that:\\
(Alm1)Low level graph $\Gamma$, $V_\Gamma= \{1,2,3\}$
 $\hspace{.3in} \Gamma \ \hspace{.3in} \ \begin{array}{ccc}
2&\frac{\hspace{.4in}}{} & 3\\
& 1 & \\
\end{array} \hspace{.3in}$

(Alm2) We have two update functions, for the family $A_1$, and for $A_2$.
   With $A_1$ we have a local function  associated to the vertex $1$,
 $f_{11}=(1,x_2,x_3)$, and for $A_2$ we have a trivial function $f_{12}=Id$,
 one associated to the second vertex $2$,  $f_{21}=Id$ for both families, and one
 for the vertex $3$, that we can find using the usual method for
 boolean functions, or the method giving in \cite{A1}.
 We have the following functions
 $F_1=\{ f_{11}=(1,x_2,x_3), f_{12}=Id\}$,   $F_2=\{f_{21}=Id\}$, and one
  for the vertex $3$, that we can find using the usual method for
  boolean functions: $F_3=\{f_{31}(x_1,x_2,x_3)=(x_1,x_2,x_2\overline{x_3})\}$.
  $f_{ij}:{\Z_2}^3\rightarrow{\Z_2}^3$, for all local function.

  (Alm3) We select an order $\alpha  = ( 1 \hspace{.1in} 2 \hspace{.1in}
  3 )$ in the set of vertices $V_\Gamma $.
  So we have two update functions
  \[f_1(x_1,x_2,x_3)=(x_1,x_2,x_2\overline {x_3})\, and\  f_2(x_1,x_2,x_3)=(1,x_2,x_2\overline
  {x_3}).\]
(Alm4) The probability $c_1=.66667$ for $f_1$, and $c_2=.3333\overline{3}$, for  $f_2$.
\end{example}

To solve the reverse engineering problem with probabilities, and sequential behavior we need to prove
the algorithm always runs for a set of data.

First, the low level graph is always possible to obtain, similarly with the sets $k_a$, and the families $A_i$, but
instead of we know that always a function with coordinate functions  acting simultaneously induce a sequential function
like our interested functions, we know that this problem is very complicated and is an open problem when we want
to preserve the number of vertex in the sequential function. For solutions of the reverse engineering problem, see \cite{A1,G, BL}.
\begin{proposition}
Let $f:k^n\rightarrow k^n$ be a function with coordinate functions $f=(f_1, \ldots,  f_n)$, and let
 $\alpha =( \alpha (1) \hspace{.1in} \cdots  \hspace{.1in} \alpha(n))$ be a permutations of the vertex of $\Gamma$.
  If in the functions $f_{\alpha (i)}$  only appears the variables $x_{\alpha (j)}$ such that $j\leq i $
  then $f=\bar{f}_{\alpha(1)} \circ \cdots \circ \bar{f}_{\alpha(n)} $, where $\bar{f}_i(x_1, \ldots, x_n)=(x_1,\ldots,x_{i-1},f_i,x_{i+1}, \ldots,
  x_n)$. So, our claim holds.
\end{proposition}
\begin{pf}
The proof is trivial. In fact, using induction and the simple case when $\alpha=(1 \ \ldots \ n)$ we have tha t
if the function $f_i$ only use the variables $x_1,\ldots ,x_i$, and the function
$f_{i-1}$ only use   $x_1,\ldots ,x_{i-1}$, then $\bar{f}_{i-1}\circ \bar{f}_{i}=(x_1,\ldots, f_{i-1},f_{i}, x_{i+1},\ldots, x_n)$.
\end{pf}
\section{Probabilistic Sequential Networks}\label{PSN1}
The following definition give us the possibility to have several
update functions acting in a sequential manner with assigned
probabilities. All these, permit us to study the dynamic of these
systems using Markov chains and other probability tools.
\begin{definition}\label{PSN}
 A Probabilistic Sequential Network (PSN) \begin{displaymath}{\mathcal
D}=(\Gamma, \{F_a\}_{a=1}^{|\Gamma|=n}, (k_a)_{a=1}^{n},
(\alpha_j)_{j=1}^m, C=\{c_1,\ldots, c_s\})\end{displaymath} consists
of:
\begin{enumerate}
\item[(1)] a finite graph
  $\Gamma=(V_\Gamma,E_\Gamma)$ with $n$ vertices;

 \item [(2)]  a family of finite  sets $(k_a|a\in V_\Gamma)$.

  \item [(3)] for each vertex $a$ of  $\Gamma$  a  set of
local functions
\[F_a=\{f_{ai}:k^n\rightarrow k^n|1\leq i \leq
\ell(i)\},\] is assigned. (i. e. there exists a bijection map $\sim
: V_\Gamma\rightarrow \{F_a |1\leq a\leq n\}$) (for definition of
local function see (\ref{SDS1}.2)).

\item [(4)] a family of $m$ permutations $\alpha  =(\alpha_{1}\hspace{.1in} \ldots \hspace{.1in} \alpha_{n}) $ in the set of vertices $V_\Gamma $.

\item [(5)]and  a set $C=\{c_1, \ldots ,c_s\}$, of
assign probabilities to $s$ update functions.
\end{enumerate}
\end{definition}
We select one function  in  each set $F_a$, that is one for each
vertices $a$ of $\Gamma$, and a permutation $\alpha $, with the
order in which the vertex $a$ is selected, so there are
$\underline{n}$ possibly different update functions
$f_i=f_{\alpha_1i_1}\circ \ldots \circ f_{\alpha_ni_n}$, where
$\underline{n}\leq n !\times\ell(1)\times \dots \times \ell (n)$.
The probabilities are assigned to the update functions, so there
exists a set $S=\{f_1,\ldots ,f_s\}$ of selected update functions
such that $c_i=p(f_i)$, $1\leq i \leq s$.
\begin{definition}
 The State Space of ${\mathcal D}$ is a weighted
digraph whose vertices are the elements of $k^{n}$ and there is an
arrow going from the vertex $u$ to the vertex  $v$  if there exists
an update function $f_i\in S$, such that $v=f_i(u)$. The probability
$p(u,v)$ of the arrow going from $u$ to $v$  is the sum of the
probabilities $c_{f_i}$ of all functions $f_i$, such that
$v=f_i(u)$, $u \hspace{.1in}
\frac{\hspace{.1in}p(u,v)\hspace{.1in}}{}> \hspace{.1in} f_i(u)=v $.
We denote the State Space by $\mathcal {S_D}$.
\end{definition} For each one update function in $S$ we have one  SDS
inside the PSN, so the State Space $\mathcal S_f$ is a subdigraph of
$\mathcal {S_D}$. When we take the whole set of update functions
generated by the data, we will say that we have the \emph{full} PSN.
It is very clear that a SDS is a  particular PSN, where we take one
local function for each vertex, and  one permutation. The dynamic of
a PSN  is described by  Markov Chains of the transition matrix
associated to the State Space.
\begin{example}\label{exam1}  Let ${\mathcal
D}=(\Gamma;F_1,F_2,F_3;\mathbf{Z_2}^{3};\alpha _1, \alpha _2;
(c_{f_i})_{i=1}^8),$ be the following PSN:\\
\emph{(1)} The graph $\Gamma$: $  \ \ \ \begin{array}{c} 1 \bullet  {\overline{\hspace{.2in}}} \bullet \ 3\\
  \diagdown   \mid \\
 2\hspace{.15in}  \bullet  \end{array}.$\\
\emph{(2)}  Let  $\textbf{x}=(x_1,x_2,x_3)\in \{0,1\}^3$. In this
paper, we  always consider the operations over the finite field
$\Z_2=\{0,1\}$, but we use additionally the following notation
$\bar{x_1}=x_1+ 1$. Then the sets of local functions from
${\Z_2}^3\rightarrow {\Z_2}^3$ are the following
\[\begin{array}{l}
F_1 =\{f_{11}(\textbf{x})=(1,x_2,x_3),f_{12}(\textbf{x})=(\bar{x_1},x_2,x_3))\}  \\
F_2 = \{f_{21}(\textbf{x})=(x_1, x_1x_2,x_3)\}   \\
F_3=\{f_{31}(\textbf{x})=(x_1,x_2,x_1x_2),f_{32}(\textbf{x})=(x_1,x_2,x_1x_2+x_3)\}
\end{array}.\]
\emph{(3)} The schedules or permutations are
$\alpha_1=( 3\hspace{.1in}2\hspace{.1in} 1); \alpha_2=( 1\hspace{.1in}2\hspace{.1in} 3) .$
  We obtain the following table of
 functions, and we select all of them for ${\mathcal D}$ because the
 probabilities given by $C$.
\[\begin{array}{ll}
f_1=f_{31}\circ f_{21} \circ f_{11} & f_2=f_{11}\circ f_{21} \circ f_{31} \\
f_3=f_{32}\circ f_{21} \circ f_{11} & f_4=f_{11}\circ f_{21} \circ
 f_{32}\\
f_5=f_{31}\circ f_{21} \circ f_{12} & f_6=f_{12}\circ f_{21} \circ
f_{31}\\
f_7=f_{32}\circ f_{21} \circ f_{12}& f_8=f_{12}\circ f_{21}\circ
f_{32} \cr
\end{array}.\]
  The  update functions are the following:
  \[\begin{array}{ll}
  f_1(\textbf{x})=(1,x_2,x_2)& f_2(\textbf{x})=(1,x_1x_2,x_1x_2) \\
 f_3(\textbf{x})=(1,x_2,x_2+x_3) & f_4(\textbf{x})=(1, x_1x_2,x_1x_2+x_3)\\
f_5(\textbf{x})=(\bar{x_1},\bar{x_1}x_2,
 (x_1+1)x_2)  & f_6(\textbf{x})=(\bar{x_1},x_1x_2,x_1x_2)\\
f_7(\textbf{x})=(\bar{x_1},(x_1+1)x_2,(x_1+1)x_2+x_3)&
f_8(\textbf{x})=(\bar{x_1},x_1x_2,x_1x_2+x_3)
 \cr
\end{array}.\]
\emph{(4)} The  probabilities assigned are the following:
$c_{f_1}=.18; c_{f_2}=.12; c_{f_3}= .18; c_{f_4}=.12; c_{f_5}=.12;
c_{f_6}=.08; c_{f_7}=.12; c_{f_8}=.08$.
\end{example}
\begin{example}\label{MA} We notice that there are several PSN that we
can construct with the same initial data of functions and
permutations, but with different set of  probabilities, that is,
subsystems of ${\mathcal D}$. For example if $S'=\{f_1, f_2,
f_3,f_4\}$, $F_1'=\{f_{11}\}$, and $D=\{ d_{f_1}=.355, d_{f_2}=.211,
d_{f_3}=.12,d_{f_4}=.314\}$, then
\[{\mathcal B}=(\Gamma;F_1',F_2,F_3;{\mathbf Z_2}^{3};\alpha _1, \alpha _2; D=\{
.355,.211, .12,.314\} ),\] is a PSN too.
\end{example}
\section{Morphisms of Probabilistic Sequential Networks}\label{secHPSS}
The definition of morphism of PSN is a natural extension of the
concept of homomorphism of SDS. In this section we prove in Theorem
\ref{C3} a strong property, that is the distribution of
probabilities of  two homomorphic PSN are enough close to prove
Theorem \ref{MT}.

Consider the following two PSN ${\mathcal
D}_1=(\Gamma,(F_a)_{a=1}^{|\Gamma|=n},(k_a)_{a=1}^{n},(\alpha
^j)_{j},C)$ and  \\ ${\mathcal D}_2
=(\Delta,(G_b)_{b=1}^{|\Delta|=m},(k_b)_{b=1}^{m},(\beta
^j)_{j},D).$ We denote by $S_i$ the set of update functions of
${\mathcal D}_i$, $i=1,2$; and the following notation for $(u,  v)
\in k^n \times k^n$, and $(w,  z) \in k^m\times k^m$,
\[c_f(u,v)=\left\{\begin{array}{cc}
p(f) & if  \ f(u)=v \\
0 & otherwise  \cr
 \end{array}\right\}, \
 d_g(w,z))=\left\{\begin{array}{cc}
p(g) & if  \ g(w)=z \\
0 & otherwise  \cr
 \end{array}\right\}\]
where $p(h)$ is the probability of the function $h$.
\begin{definition}{\emph{(Morphisms of PSN)}}\label{homPSS} A
morphism $h:{\mathcal D}_1\rightarrow {\mathcal D}_2$  consist of:
\begin{enumerate}
\item [(1)]  A graph morphism  $\phi:\Delta \rightarrow \Gamma $, and a family of maps in the category \textbf{Set},
$(\widehat{\phi}_{b}:k_{\phi(b)}\rightarrow k_{b}\forall b\in \Delta
),$ that  induces the adjoint function $h_\phi$, see (\ref{defh}).
 \item [(2)]  The induced adjoint map $h_\phi: k^{n}\rightarrow k^{m}$ holds that
 for all update functions $f$ in $S_1$ there exists an
update function $g\in S_2$ such that $h$ is a SDS-morphism from
$(\Gamma,(f:k^n\rightarrow k^n),\alpha_j)$ to
$(\Delta,(g:k^m\rightarrow k^m),\beta _j)$. That is, the diagrams
\ref{eq1}, \ref{eq2}, and \ref{eq3} commute for all  $f$ and its
selected $g$.
\begin{equation}\label{eq4}
h_\phi \circ f_{\alpha_1}\circ \cdots \circ f_{\alpha_n}
=g_{\beta_1}\circ \cdots \circ g_{\beta_m} \circ h_\phi
\end{equation}
 \end{enumerate}
\end{definition}
The second condition induces a  map $\mu $ from $S_1$ to $S_2$, that
is $\mu (f)=g$  if the selected function for $f$ is $g$.
 We say that a morphism
$h$ from ${\mathcal D}_1$ to ${\mathcal D}_2$ is a
\textbf{PSN-isomorphism} if $\phi$,  $h_\phi$, and $\mu$ are
bijective functions, and $d(h_\phi(u),h_\phi(g(u))= c(u,f(u))$ for
all $u$, in $k^n$, and all $f\in S_1$, and all $g\in S_2$. We denote
it by ${\mathcal D}_1\cong {\mathcal D}_2$. \vspace{.06in}

\hspace{2in}\textsc{Some theorems}
\begin{theorem}\label{C3} The morphism $h:\mathcal
D_1\rightarrow\mathcal D_2$ induces the following probabilistic
condition:

 For a fixed real number $0\le \epsilon<1$,  the map
$h_\phi$ satisfies the following:
\begin{equation}\label{epsilon}
 max_{u,v}|c_{f}(u,v)-d_{g}(h_\phi(u),h_\phi(v))|\le \epsilon
 \end{equation}
 for all  $f$ in $S_1$, and its selected
$g$ in $S_2$, and all $(u, v)\in k^n \times k^n$.
\end{theorem}
\begin{pf} Suppose $\phi,$ and $h_\phi$ satisfy  the Definition
\ref{homPSS}; and \[|c_{f}(u,v)-d_{g}(h_\phi(u),h_\phi(v))|\ge 1\]
for some $(u,v)\in k^n\times k^n$. Then we have one of the following
cases
\begin{itemize}
\item [1.] $c_{f}(u,v)=1$ and
$d_{g}(h_\phi(u),h_\phi(v))=0$. It is impossible by  condition (2)
in  definition \ref{homPSS}. In fact, if we have an arrow going from
$u$ to $v=f(u)$, then  there exists an arrow going from $h_\phi(u)$
to $h_\phi(v)=g(h_\phi(u))$ by diagram \ref{eq4}, and the
probability $d_{g}(h_\phi(u),h_\phi(v))\ne 0$.
\item[2.] $c_{f}(u,v)=0$, and $d_{g}(h_\phi(u),h_\phi(v))=1$.
 It is impossible because at least there exists  one element
$v_1\in k^n$, such that $f(u)=v_1\in k^{n}$ and $c_{f}(u,v_1)\ne 0$,
then $d_{g}(h_\phi(u),h_\phi(v_1))\ne 0$ too. Since the sum of
probabilities of all arrow going up from $h_\phi(u)$ is equal $1$,
then $d_{g}(h_\phi(u),h_\phi(v))<1$, and our claim holds.
\end{itemize}
Therefore the  condition holds, and always $\epsilon$ exists.
\end{pf} In the next theorem we will use the following notation:
\begin{itemize}
\item [(1)]  $S_\phi=\mu (S_1)$.
\item [(2)] $\delta^{t}=\sum_{g \not \in
S_\phi}d_{g^{t}}$, where $g^{t}=g\circ g \circ \cdots \circ g$, $t$
times.
\item [(3)] $p_t(u,v)=\sum_{f^{t}}c_{f^{t}}(u,v)$, and $p_t(h_\phi(u),h_\phi(v))=\sum_{g^{t}}d_{d^{t}}(h_\phi(u),h_\phi(v))$
\item[(4)] $T_i$ denotes the transition matrix of the PSN ${\mathcal
D}_i$, for $i=1,2$, and \\ $p_t(u,v)=({T_i}^t)_{(u,v)}$.
\end{itemize}
 \begin{theorem}\label{MT}
If $h:{\mathcal D}_1\longrightarrow{\mathcal D}_2$ is a morphism of
probabilistic sequential networks, then:
\[\lim
_{m\rightarrow\infty}|({T_1})^{m}_{u,v}-({T_2})^{m}_{h_\phi(u),h_\phi(v)}|=0,\]
 for all   $(u,v)\in k^n \times k^n$. That is, the equilibrium state
 of both systems are equals.
 \end{theorem}
\begin{pf}
The condition giving by Theorem \ref{C3} asserts that, there exists
a fixed real number $0\le \epsilon<1$,  such that the map $h_\phi$
satisfies the following:
\[
 {\max}_{u,v}|c_{f}(u,v)-d_{g}(\phi(u),\phi(v))|\le \epsilon
 \]
 for all  $f$ in $S_1$, and its selected
$g$ in $S_2$, and  all $(u, v)\in k^n \times k^n$.

If there is a function $f$ going from $u$ to $v=f(u)$  in $k ^{n}$,
then there exists a function $g$ going from $h_\phi(u)$ to
$h_\phi(v)$, such that  $g(h_\phi(u))=h_\phi(f(u))$. Now, for $m=
2$,
 and by the Chapman-Kolmogorov equation \cite{Ste},  the
 following  computation is valid
\[|c_{f^2}(u,f^2(u))-d_{g^2}(h_\phi(u),g^2(h_\phi(u)))|=\]
\[|c_{f}(u,f(u))c_{f}(f(u),f^2(u))-d_{g}(h_\phi(u),g(h_\phi(u)))d_{g}(g(h_\phi(u)),g^2(h_\phi(u)))|=\]
\[|c_{f}(u,f(u))c_{f}(f(u),f^2(u))-d_{g}(h_\phi(u),h_\phi(f(u)))d_{g}(h_\phi(f(u)),h_\phi(f^2(u)))|\leq\]
\[\leq
|c_{f}(f(u),f^2(u))|\epsilon+|d_{g}(h_\phi(u),h_\phi(f(u)))|\epsilon\leq
2\epsilon,\]  by condition \ref{epsilon}. We just proved that
$|c_{f^2}(u,f^2(u))-d_{g^2}(h_\phi(u),g^2(h_\phi(u)))|\leq 2
\epsilon$.
 Using  mathematical induction  over $m$, we conclude that, for
 all natural number $m\geq 1$
\begin{equation}
{\max}_{u,f^{m}(u)}|(c_{f^m}(u,f^m(u))-d_{g^m}(h_\phi(u),g^m(h_\phi(u)))|\leq
m\epsilon . \label{eq:1}
\end{equation}
For $m=2$,  this result implies that
\[|p_2(u,v)-p_2(h_\phi(u),h_\phi(v))|\leq 2k \epsilon +
\delta^{2},\] where $k$ is the maximum number of functions ${f}^{2}$
going from one state to another in $k^{n}$. The sum $\delta^{2}$ is
taking over  the functions $g$ that are going from $h_\phi(u)$ to
$h_\phi(v)$ and do not  have a  function $f$ in $S_1$  associated to
the function $g$. So, the sum is not over  all the functions in
$S_2$, and  we have $\delta^{2} < 1$, and maybe $\delta^{2} =0$, see
\cite{Ste}. Then the above condition  implies that:
\begin{equation}
   \hspace{.3in} {\max}_{(u,v)\in k^n\times k^n}|p_2(u,v)-p_2(h_\phi(u),h_\phi(v))|\leq
2k\epsilon + \delta ^{2}  \label{eq:2}
\end{equation}
 Using induction, we conclude that
\begin{equation}
  \hspace{.3in} {\max}_{(u,v)\in k^n\times k^n}|p_m(u,v)-p_m(h_\phi(u),h_\phi(v))|\leq
mk\epsilon + \delta ^{m}  \label{eq:3}
 \end{equation}   for all  $m\in \mathbf N $, the natural numbers.

So, for all real number $0<\epsilon ' < 1$ there exists $m\in\mathbf
N$, such that,
\[|p_{m'}(u,v)-p_{m'}(h_\phi(u),h_\phi(v))|< \epsilon ',\]
 for all natural number  $m'>m$, and for all
possible $u, \ v \in k^n$.

 In fact, we have $\epsilon ^{m'}\ll \epsilon
 ^{m}$, and this implies
 $(m'k)\epsilon ^{m'}<(mk) \epsilon ^{m}$. Similarly
  $\delta ^{m'}< \delta ^{m}$, so  selecting $m$ such that
   $(mk)\epsilon ^{m}+ \delta ^{m}< \epsilon '$, we obtain
\[|p_{m'}(u,v)-p_{m'}(h_\phi(u),h_\phi(v))|\leq (m'k') \epsilon ^{m'}+ \delta ^{m'}< (mk)\epsilon^{m}+ \delta ^{m}< \epsilon ',\]
where $k'$ is the maximum number of functions going from one state
to another in the state space of the power $m'$ of the functions. We
can observe that because the state space $\mathcal S_\mathcal D$ is
finite, $k'\leq k$. Therefore
\[\lim
_{m\rightarrow\infty}|p_{m'}(u,v)-p_{m'}(h_\phi(u),h_\phi(v))|=0,\]
 for all possible  $(u,v)\in k^n\times k^n$, and
 the theorem holds.
\end{pf}
\begin{corollary}
Two probabilistic sequential network are homomorphic if they have
the same probabilistic equilibrium.
\end{corollary}
\begin{pf}
It is obvious using the proof of the theorem, because if they have
the same probabilistic equilibrium, then the two Time Discrete
Markov Chains have the same size. On the other hand, $\delta$ is
almost $0$ and there exists a morphism going from one PBN to the
other.
\end{pf}

\hspace{-.1in}\textsc{Special morphisms.} Let ${\mathcal
D}=(\Gamma,(F_i)_{i=1}^n,(\alpha ^j)_{j\in{J}},C)$ be a PSN.

\hspace{-.1in}\textsc{Identity morphism.}
  The functions $\phi=id_\Gamma,$ $h_\phi=id_{k^n}$  and
$\mu =id_{S}$, define  the \emph{identity morphism} ${\mathcal
I}:\mathcal D\rightarrow \mathcal D$, and it is a trivial example of
a PSN-isomorphism.

\hspace{-.1in}\textsc{Monomorphism} A morphism $h$ of PSN is a
\emph{ monomorphism} if $\phi$ is surjective and $h_{\phi}$ is
injective.

\hspace{-.1in}\textsc{Epimorphism} A morphism is an
\emph{epimorphism} if $\phi$ is injective and $h_{\phi}$ is
surjective.

\hspace{-.1in}\textsc{Remark} If the morphism $h$ is either a
monomorphism or an epimorphism, then the function $\mu$ is not
necessary injective, neither surjective.
\section{The category  \textbf{PSN}. Simulation in the category \textbf{PSN}, and Examples}\label{CAT}
In this section, we  prove that the PSN with the \emph{morphisms}
form a category, that we denote by \textbf{PSN}. For unknown
definitions, and results in Categories see the famous and old book
of S. MacLane: Categories for the Working Mathematicians \cite{ML}.
\begin{theorem}\label{COM} Let
$ h_1=(\phi _1, h_{\phi _1}):{\mathcal D}_1\rightarrow {\mathcal
D}_2$ and   $h_2=(\phi _2, h_{\phi _2}):{\mathcal D}_2\rightarrow
{\mathcal D}_3$ be two morphisms of PSN. Then the composition
$h=(\phi , h_{\phi })=(\phi _2, h_{\phi _2}) \circ (\phi _1, h_{\phi
_1})=h_2\circ h_1: {\mathcal D}_1\rightarrow {\mathcal D}_3$  is
defined as follows: $h=(\phi , h_{\phi })=(\phi_1 \circ \phi _2,
h_{\phi _2} \circ h_{\phi _1})$ is  a morphism of PSN. The function
$\mu_h=\mu_{h_2}\circ \mu_{h_1}.$
\end{theorem}
\begin{pf}
The composite function $\phi=\phi_1 \circ \phi_2 $ of two graph
morphisms is again a graph morphism. The composite function
$h_\phi=h_{\phi _2} \circ h_{\phi _1}$ is again a digraph morphism
which satisfies the conditions in Definition \ref{homPSS}, by
Proposition and Definition 2.7 in \cite{LP2}.  So, $h=(\phi,
h_\phi)$ is again a morphism. of PSN.
\end{pf}
\begin{theorem}
The Probability Sequential Networks together with the homomorphisms
of PSN form the category \textbf{PSN}.
\end{theorem}
\begin{pf} Easily  follows from Theorem \ref{COM}.
\end{pf}
\begin{theorem}\label{FSDS} The SDS together with the morphisms defined in \cite{LP2} form a full subcategory of the category
\textbf{PSN}.
\end{theorem}
\begin{pf}
It is trivial.
\end{pf}
 \subsection{Simulation and  examples}\label{EHOM}
In this section we give several examples of morphisms, and
simulations. In the second example we show how the Definition
\ref{homPSS} is verified under the supposition that a function
$\phi$ is defined. So, we have two examples in (\ref{EHOM}.2), one
with $\phi$ the natural inclusion, and the second  with $\phi$  a
surjective map. The third, and the fourth examples are morphisms
that represent simulation of $\mathcal G$ by $\mathcal F$. We begin
this section with the definitions of Simulation and sub-PSN.

\hspace{-.15in}\textsc{Definition of Simulation in the category
\textbf{ PSN}}.
 The probabilistic sequential network $\mathcal G$ is simulated by $\mathcal F$ if
there exists a  monomorphism $h: {\mathcal F}\rightarrow {\mathcal
G}$ or an epimorphism $h':{\mathcal G}\rightarrow {\mathcal F}$.

\hspace{-.15in}\textsc{Sub Probabilistic Sequential
Network}\label{subPSS}
 We say that a PSN ${\mathcal G}$ is a
sub Probabilistic Sequential Network of ${\mathcal F}$ if there
exists a monomorphism from ${\mathcal G}$ to ${\mathcal F}$. If the
map $\mu$ is  not a bijection, then we say that it is a proper
sub-PSN.

\subsection{Examples}\label{examples}\hspace{2.2in}

 \textbf{(\ref{examples}.1)} In the examples \ref{exam1}, and \ref{MA} we
define  two PSN $\mathcal D$ and $\mathcal B$. The functions $\phi
=Id_\Gamma$, $h_\phi=Id_{{k}^{n}}$, and $\mu$ the natural inclusion
from $S_1$ to $S_2$ define the inclusion $\iota_\mu: {\mathcal
B}\rightarrow {\mathcal D}$. It is clear that the inclusion is a
monomorphism, so  $\mathcal D$ is simulated by ${{\mathcal B}}$.

\textbf{(\ref{examples}.2)} We now construct a monomorphism $h:
\mathcal F\rightarrow \mathcal G$, with the properties that $\phi$
 is surjective and the function $h_\phi$ is  injective.
   The PSN ${\mathcal F}=(\Gamma ,(F_i
)_3,\beta ,C )$ has the support graph $\Gamma$ with three vertices,
and the PSN ${\mathcal
 G}=(\Delta ,(G_i )_4,\alpha,D)$ has the support graph $\Delta $ with four
vertices
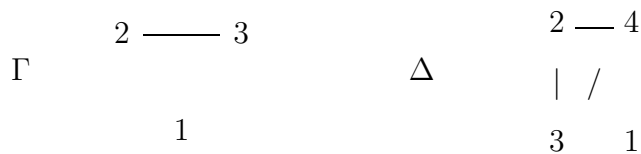
\[\hspace{.3in} \Gamma \ \hspace{.3in} \ \begin{array}{ccc}
2&\frac{\hspace{.4in}}{}  &  3\\
& 1 &  \\
\end{array} \hspace{.3in}\ \  \hspace{.4in} \Delta \ \ \hspace{.3in} \ \ \
\begin{array}{ccc}
\vspace{-.2in}
 2 &{\overline{\hspace{.2in}}} &  4\\
 \vspace{-.2in}
 \mid & / &  \\
 3& \hspace{.2in} &  1
\end{array}  \hspace{.2in}   \    \]
The morphism $h:{\mathcal F} \rightarrow {\mathcal G}$, has the
contravariant  graph morphism $\phi : \Delta\rightarrow \Gamma$,
defined by the arrows of graphs, as follows $\phi(1)=1$,
$\phi(2)=\phi(3)=2,$ and $\phi(4)=3 $, so it is a surjective map.
The family of functions $\hat{\phi}_i:{k}_{\phi (i)}\rightarrow
 k_{(i)}$,  $\hat{\phi}_1(x_1)=x_1$; $\hat{\phi}_2(x_2)=x_2;\
 \hat{\phi}_3(x_2)=x_2;\  \hat{\phi}_4(x_4)=x_4$, are injective functions.
 The sets $k_a=\Z_2$, for all vertices $a$ in $\Delta$, and $\Gamma$. The adjoint function is
 $h_\phi:{\mathbf  Z_2}^3\rightarrow {\mathbf  Z_2}^4,$  defined by
  \[  h_\phi(x_1,x_2,x_3)=(\hat{\phi}_1(x_1),\hat{\phi}_2(x_2),\hat{\phi}_3(x_2),\hat{\phi}_4(x_4))=(x_1,x_2,x_2,x_3).\]
 Then, the first condition in the definition \ref{homPSS} holds.

The  PSN  ${\mathcal F}=(\ \Gamma ;\ (F_i)_3 ;\ \beta  ;\ C
),$ is defined with the following data.\\
The set of  functions $ F_1=\{ \ f_{11}, \ f_{12})\},$
 $F_2=\{f_{21}\},$ and $
F_3=\{f_{31}\},$ where
\[f_{11}=Id,\ f_{12}(x_1,x_2,x_3)=(1,x_2,x_3), \
f_{21}=Id, \]
\[f_{31}(x_1,x_2,x_3)=(x_1,x_2,x_2\overline{x_3}).\]
A permutation  $\beta =(\ 1 \ 2 \ 3\ );$ and the probabilities
$C=\{c_{f_1}=.5168, c_{f_2}=.4832\} $. So, we are taking all the
update functions   $S =\{f_1,f_2\}$;
 \[f_1= f_{11}\circ f_{21}\circ f_{31},\
f_1=(x_1,x_2,x_3)= (x_1,x_2,x_2\overline{x_3});\]
\[\textrm{and} \ \ \ f_2= f_{12}\circ f_{21}\circ f_{31}, \
f_2(x_1,x_2,x_3)= (1,x_2,x_2\overline{x_3}).\]

On the other hand, the PSN ${\mathcal G}=(\Delta ;(G_i)_4;\alpha; D )$ has the following data.\\
The families of functions: $G_1=\{g_{11}, g_{12}\}$; $G_2=\{g_{21},
g_{22}\}$, $G_3=\{g_{31},g_{32}\}$; and $G_4=\{g_4\}$, where
\[\begin{array}{ll}
g_{11}(x_1,x_2,x_3,x_4)&= (1,x_2,x_3,x_4) \\
g_{21}(x_1,x_2,x_3,x_4)&= (x_1,1,x_3,x_4)  \\
g_{31}(x_1,x_2,x_3,x_4)&= (x_1,x_2,x_1x_2,x_4) \\
g_{41}(x_1,x_2,x_3,x_4)&=(x_1,x_2,x_3,x_2\overline{x_4})\\
\end{array}\ \begin{array} {ll}
g_{12}= Id=g_{22} &  \\
g_{32}(x_1,x_2,x_3,x_4)&= (x_1,x_2,x_2,x_4)\\
\end{array}.\]
One permutation or schedule  $\alpha=(1\hspace{.1in}2\hspace{.1in}
3\hspace{.1in}4)$. The assigned probabilities $d_{g_5}=
.00252,\ d_{g_6}=.08321,\  d_{g_7}=.51428,\  d_{g_8}=.39999$ whose
determine the set of update functions $X=\{g_5,g_6,g_7,g_8\}$:
therefore the all update functions are  the following
\[\begin{array}{lll}
g_1=g_{11} \circ g_{21} \circ g_{31 }\circ g_{41} ,& g_2= g_{12}
\circ g_{21} \circ g_{32}\circ g_{41} &
g_3=g_{12} \circ g_{21} \circ g_{31}\circ g_{41},\\
 g_4=g_{11}\circ g_{21}\circ g_{32}\circ g_{41}& g_5=g_{12} \circ g_{22} \circ g_{31}\circ g_{41},& g_6=g_{11} \circ g_{22} \circ g_{31} \circ g_{41}\\
g_7=g_{12} \circ g_{22} \circ g_{32}\circ g_{41},& g_8=g_{11}\circ
g_{22}\circ g_{32}\circ g_{41} &\cr
\end{array}.\]
The selected functions are
\[\begin{array}{ll}
g_5(x_1,x_2,x_3,x_4)=(x_1,x_2,x_1x_2,x_2\overline{x_4}),&
g_6(x_1,x_2,x_3,x_4)=(1,x_2,x_1x_2,x_2\overline{x_4})\\
g_7(x_1,x_2,x_3,x_4)=(x_1,x_2,x_2,x_2\overline{x_4}),&
g_8(x_1,x_2,x_3,x_4)=( 1,x_2,x_2,x_2\overline{x_4})\\
\end{array}.\]
  We claim that $h: \mathcal F \rightarrow \mathcal G$ is a morphism.
  It is trivial that the following diagrams commute.
\[\begin{array}{lcl}
{\mathbf  Z_2}^3 &\overrightarrow{\hspace{.2in}_{f_1}\hspace{.2in}}& {\mathbf  Z_2}^3 \\
h_\phi \downarrow &  & \downarrow h_\phi \\
 {\mathbf  Z_2}^4 & \overrightarrow{\hspace{.2in}_{g_7}\hspace{.2in}} & {\mathbf  Z_2}^4 \cr
\end{array}, \ \textrm{and} \  \ \begin{array}{lcl}
{\mathbf  Z_2}^3 & \overrightarrow{\hspace{.2in}_{f_2}\hspace{.2in}}& {\mathbf  Z_2}^3 \\
h_\phi \downarrow &  & \downarrow h_\phi \\
 {\mathbf  Z_2}^4 &\overrightarrow{\hspace{.2in}_{g_8}\hspace{.2in}} & {\mathbf  Z_2}^4 \cr
\end{array}.\]
In fact, $(h_\phi \circ f_1)(x_1, x_2,x_3)
=(x_1,x_2,x_2,x_2\overline{x_3})=(g_7\circ h_\phi)(x_1, x_2,x_3),$
 on the other hand $(h_\phi \circ
f_2)(x_1,x_2,x_3)=(1,x_2,x_2,x_2\overline{x_3})=(g_8\circ
h_\phi)(x_1,x_2,x_3)$ so, the property holds.  We verify the second
property in the definition of morphism  for the compositions $f_1$
and $g_7$, and also with the compositions $f_2$ and $g_8$. That is,
we check the sequence of local functions too.
\[\begin{array}{ccccccc}
 {\mathbf  Z_2}^3 &\overrightarrow{\hspace{.2in}_{f_{31}}\hspace{.2in}} & {\mathbf  Z_2}^3 &\overrightarrow{\hspace{.2in}_{f_{21}}\hspace{.2in}}&
 {\mathbf  Z_2}^3 &\overrightarrow{\hspace{.2in}_{f_{11}}\hspace{.2in}} & {\mathbf  Z_2}^3  \\
  h_\phi\downarrow & & h_\phi \downarrow& & h_\phi\downarrow &
&h_\phi\downarrow \\
 {\mathbf  Z_2}^4 &\overrightarrow{\hspace{.2in}_{g_{41}}\hspace{.2in}} & {\mathbf  Z_2}^4 &\overrightarrow{\hspace{.2in}_{g_{22}\circ g_{32}}\hspace{.2in}}
& {\mathbf  Z_2}^4
&\overrightarrow{\hspace{.2in}_{g_{12}}\hspace{.2in}} & {\mathbf
Z_2}^3
\end{array}\]

\[\begin{array}{ccccccc}
 {\mathbf  Z_2}^3 &\overrightarrow{\hspace{.2in}_{f_{31}}\hspace{.2in}} & {\mathbf  Z_2}^3 &\overrightarrow{\hspace{.2in}_{f_{21}}\hspace{.2in}}& {\mathbf  Z_2}^3
&\overrightarrow{\hspace{.2in}_{f_{12}}\hspace{.2in}}& {\mathbf  Z_2}^3  \\
  h_\phi\downarrow & & h_\phi \downarrow& & h_\phi\downarrow &
&h_\phi\downarrow \\
 {\mathbf  Z_2}^4 &\overrightarrow{\hspace{.2in}_{g_{41}}\hspace{.2in}} & {\mathbf  Z_2}^4 &\overrightarrow{\hspace{.2in}_{g_{22}\circ g_{32}}\hspace{.2in}}
& {\mathbf  Z_2}^4
&\overrightarrow{\hspace{.2in}_{g_{11}}\hspace{.2in}} & {\mathbf
Z_2}^3
\end{array}\]
$(h_\phi \circ
f_{31})(x_1,x_2,x_3)=(x_1,x_2,x_2,x_2\overline{x_3})=(g_{41}\circ
h_\phi)(x_1,x_2,x_3)$, \\
$(h_\phi \circ f_{21})(x_1,x_2,x_3)=(x_1,x_2,x_2,x_3)=((g_{22} \circ
g_{32}) \circ
h_\phi)(x_1,x_2,x_3)$,\\
$(h_\phi \circ f_{11})(x_1,x_2,x_3)=(x_1,x_2,x_2,x_3)=(g_{12}\circ
h_\phi)(x_1,x_2,x_3)$,\\
 $(h_\phi \circ f_{12})(x_1,x_2,x_3)=(1,x_2,x_2,x_3)=(g_{11}\circ
h_\phi)(x_1,x_2,x_3)$. \\
 Then our claim holds.

 \textbf{(\ref{examples}.3)} We can construct
an epimorphism $h':\mathcal G\rightarrow \mathcal F$,  that is, the
function $\phi$ is injective and the function $h'_{\phi}$ is
surjective. We use $\phi':\Gamma\rightarrow\Delta$, defined as
follow $\phi'(i)=i+1$, for all $i\in V_\Gamma$. Therefore
$\hat{\phi'}_i:{k}_{\phi'(i)}\rightarrow k_{(i)}$,
$\hat{\phi'}_i:\Z_2\rightarrow  \Z _2$, for all $i\in V_\Gamma $,
and should be satisfies
$<h_\phi(x),i>:=\hat{\phi}_b(<x,\phi(i)>)=\hat{\phi}_b(x_{\phi(i)})$.
So, the adjoint function is
$h'_\phi(x_1,x_2,x_3,x_4)=(\hat{\phi'}_1(x_1),\hat{\phi'}_2(x_3),\hat{\phi'}_3(x_3))=(x_1,x_2,
x_4)$ and  satisfies the following commutative diagrams
\[\hspace{.4in}\begin{array}{lll}
{\Z_2}^4 & \overrightarrow{\hspace{.2in}_{g_{5}}\hspace{.2in}} &{\Z_2}^4 \\
 \downarrow h'_{\phi}&   & \downarrow h'_{\phi}   \\
{\Z_2}^3& \overrightarrow{\hspace{.2in}_{f_1}\hspace{.2in}} &{\Z_2}^3\\
\end{array}, \hspace{.4in}  \begin{array}{lll}
{\Z_2}^4 & \overrightarrow{\hspace{.2in}_{g_{7}}\hspace{.2in}} &{\Z_2}^4 \\
\downarrow h'_{\phi}&   & \downarrow h'_{\phi}   \\
{\Z_2}^3& \overrightarrow{\hspace{.2in}_{f_1}\hspace{.2in}} &{\Z_2}^3\\
\end{array}, \]
 \[\hspace{.4in}\begin{array}{lll}
{\Z_2}^4 &\overrightarrow{\hspace{.2in}_{g_{6}}\hspace{.2in}}  &{\Z_2}^4 \\
 \downarrow h'_{\phi}&   & \downarrow h'_{\phi}   \\
 {\Z_2}^3 & \overrightarrow{\hspace{.2in}_{f_{2}}\hspace{.2in}}&{\Z_2}^3 \\
 \end{array}, \ \hspace{.4in}  \begin{array}{lll}
{\Z_2}^4 & \overrightarrow{\hspace{.2in}_{g_{8}}\hspace{.2in}} &{\Z_2}^4 \\
\downarrow h'_{\phi}&   & \downarrow h'_{\phi}   \\
{\Z_2}^3& \overrightarrow{\hspace{.2in}_{f_{2}}\hspace{.2in}} &{\Z_2}^3\\
\end{array}.\]
These implies that $\mu(g_5)=\mu(g_7)=f_1$, and
$\mu(g_6)=\mu(g_8)=f_2$. In fact,
\[(h'_\phi \circ g_5)(x_1,
x_2,x_3,x_4)=(x_1,x_2,x_2\overline{x_4})=(f_1\circ
h'_\phi)(x_1,x_2,x_3,x_4),\]
\[(h'_\phi \circ g_6)(x_1,x_2,x_3,x_4)=(1,x_2,x_2\bar{x_4})=(f_2\circ
h'_\phi)(x_1,x_2,x_3,x_4), \]
\[(h'_\phi \circ g_7)(x_1,
x_2,x_3,x_4)=(x_1,x_2,x_2\overline{x_4})=(f_1\circ
h'_\phi)(x_1,x_2,x_3,x_4),\]
\[(h'_\phi \circ g_8)(x_1,
x_2,x_3,x_4)=(1,x_2,x_2\bar{x_4})=(f_2\circ
h'_\phi)(x_1,x_2,x_3,x_4).\]

Checking the compositions of local functions $g_5=g_{12} \circ
g_{22} \circ g_{31}\circ g_{41}$, and \\$f_1= f_{11}\circ
f_{21}\circ f_{31}$, we have that the following diagrams commute
\[\begin{array}{ccccccccc}
  {\Z_2}^4& \frac{\hspace{.1in} g_{12} \hspace{.1in}}{}> & {\Z_2}^4& \frac{\hspace{.1in} g_{22} \hspace{.1in}}{} > &{\Z_2}^4 &
  \frac{\hspace{.1in} g_{31} \hspace{.1in}}{} >& {\Z_2}^4 & \frac{\hspace{.1in} g_{41} \hspace{.1in}}{} >& {\Z_2}^4 \\
 h'_{\phi}\downarrow  & & \downarrow  h'_{\phi}  & &\downarrow h'_{\phi} & &\downarrow h'_{\phi}& &\downarrow h'_{\phi} \\
 {\Z_2}^3 &\frac{\hspace{.1in} f_{11} \hspace{.1in}}{} >& {\Z_2}^3 &\frac{\hspace{.1in} Id \hspace{.1in}}{} >& {\Z_2}^3  &
 \frac{\hspace{.1in} f_{21} \hspace{.1in}}{} > &
 {\Z_2}^3 &\frac{\hspace{.1in} f_{31} \hspace{.1in}}{}> &{\Z}^2\\
\end{array}.\]

By the data we only need to check the following compositions \\
$h'_\phi(g_{31}(x_1,x_2,x_3,x_4))=(x_1,
x_2,x_4) =f_{21}(h'_\phi(x_1,x_2,x_3,x_4)),$\\
$h'_\phi(g_{41}(x_1,x_2,x_3,x_4))=(x_1, x_2,x_2\bar{x_4})
=f_{31}(h'_\phi(x_1,x_2,x_3,x_4)).$ Similarly, we can prove that the
other functions hold the property.
\section{ Equivalent  Probabilistic Sequential Networks}\label{epsn}
\begin{definition}{\emph{(Equivalent PSN)}}\label{equiv} If the  morphism  $h:\mathcal D_1 \rightarrow \mathcal D_2$  satisfies
that $\phi$,  $h_\phi$  and $\mu$ are bijective functions, but the
probabilities are not necessary equals, we say that ${\mathcal
D}_1$, and ${\mathcal D}_2$ are equivalent PSN. We write ${\mathcal
D}_1 \simeq  {\mathcal D}_2$.
\end{definition}
So, ${\mathcal D}_1$, and ${\mathcal D}_2$ are equivalents if there
exist $(\phi,h_\phi,\mu)$, and $(\phi^{-1}, h_\phi^{-1},\mu ^{-1})$,
such that for all update functions $f\in {\mathcal D}_1$ and its
selected function $g\in {\mathcal D}_2$, the condition
$f=h_\phi^{-1}\circ g \circ h_\phi$ holds . It is clear that this
relation is an equivalence relation in the set of PSN.
\begin{proposition} If ${\mathcal D}_1 \simeq  {\mathcal D}_2$, then the
transition matrices $T_1$ and $T_2$ satisfy: $(T_1^m)_{(u,v)}\ne 0$,
if and only if $(T_2^m)_{(h_\phi{(u)},h_\phi(v))}\ne 0$, for all
$m\in \N$, $(u,v)\in k^n\times k^n$.
\end{proposition}
\begin{pf} It is obvious.
\end{pf}
\section*{Acknowledgment} This work  was supported by the Partnership M. D.
Anderson Cancer Center, University of Texas, and the Medical
Sciences Cancer Center, University of Puerto Rico, by  the Program
AABRE of Rio Piedras Campus,  and by the SCORE Program of NIH, Rio
Piedras Campus.

\end{document}